\documentclass[]{interact}

\usepackage{epstopdf}

\usepackage[doublespacing]{setspace}
\usepackage[numbers,sort&compress]{natbib}
\bibpunct[, ]{[}{]}{,}{n}{,}{,}
\usepackage{csquotes}
\usepackage{epsfig}
\usepackage{lscape}
\usepackage[misc,geometry]{ifsym}
\usepackage{color}
\usepackage{booktabs}
\usepackage{multirow}
\usepackage{siunitx}
\usepackage{amsmath,amssymb,dsfont}
\usepackage{listings}
\usepackage{rotating}
\usepackage{multicol}
\usepackage[english]{babel}
\usepackage{lineno}
\usepackage{enumerate}
\usepackage[justification=centering]{caption}
\usepackage{subcaption}
\newcounter{magicrownumbers}
\newcommand\rownumber{\stepcounter{magicrownumbers}\arabic{magicrownumbers}}
\theoremstyle{definition}
\usepackage{csquotes}
\usepackage{epsfig}
\usepackage{lscape}
\usepackage{amsfonts}
\usepackage[misc,geometry]{ifsym}
\usepackage{color}
\usepackage{amsmath,amssymb,dsfont}
\usepackage{listings}
\usepackage{rotating}
\usepackage{multicol,multirow}
\usepackage[english]{babel}
\usepackage{lineno}
\usepackage{enumerate}
\usepackage{caption}
\usepackage{afterpage}
\usepackage{empheq}
\newcommand{\bea}{\begin{eqnarray}}
\newcommand{\eea}{\end{eqnarray}}
\newcommand{\nn}{\nonumber}
\newcommand{\bee}{\begin{eqnarray*}}
\newcommand{\eee}{\end{eqnarray*}}

\newcommand{\bt}{\begin{tabbing}}
\newcommand{\et}{\end{tabbing}}
\newcommand{\btb}{\begin{tabular}}
\newcommand{\etb}{\end{tabular}}

\newcommand{\bc}{\begin{center}}
\newcommand{\ec}{\end{center}}

\newtheorem{t1}{Theorem}
\newcommand{\lb}{\label}

\newtheorem{n1}{Note}
\newtheorem{d1}{Definition}
\newtheorem{l1}{Lemma}

\newtheorem{alg1}{Algorithm}
\begin{document}


\title{A Newton-Type Proximal Gradient Method for Nonlinear Multi-objective Optimization Problems}

\author{
\name{Md Abu Talhamainuddin Ansary\textsuperscript{a}\thanks{Md Abu Talhamainuddin Ansary Email: md.abutalha2009@gmail.com}}
\affil{\textsuperscript{a}Department of Economic Sciences,\\ Indian Institute of
              Technology Kanpur,\\ India-208 016}
}

\maketitle

\begin{abstract}
In this paper, a globally convergent Newton-type proximal gradient method is developed for composite multi-objective optimization problems where each objective function can be represented as the sum of a smooth function and a nonsmooth function. The proposed method deals with unconstrained convex multi-objective optimization problems. This method is free from any kind of priori chosen parameters or ordering information of objective functions. At every iteration of the proposed method, a subproblem is solved to find a suitable descent direction. The subproblem uses a quadratic approximation of each smooth function. An Armijo type line search is conducted to find a suitable step length. A sequence is generated using the descent direction and step length. The Global convergence of this method is justified under some mild assumptions. The proposed method is verified and compared with some existing methods using a set of problems. 
\end{abstract}

\begin{amscode}
90C25;  90C29; 49M37; 65K10
 \end{amscode}
\section{Introduction}
\label{intro}
In a multi-objective optimization problem, several objective functions are minimized simultaneously. If any feasible solution minimizes all objective functions, then this is called an ideal solution. But quite often decrease of one objective function causes an increase in other objective functions. So the concept of optimality is replaced by efficiency. Classical methods of solving multi-objective optimization problems are scalarization methods (see \cite{kd0,erg1,eichfelder1}), which reduce the original problem to a single objective optimization problem using a set of priori chosen parameters. These methods are user dependent and often fail to generate Pareto front. Heuristic methods like evolutionary algorithms (see \cite{deb2000fast,deb2013evolutionary,coello1,seada2018non,deb2017recent}), are often used to find the approximate Pareto front but cannot guarantee any convergence property.\\

Recently many researchers have developed new techniques for nonlinear multi-objective optimization problems that do not involve any priori chosen parameters or ordering information of objective functions. In addition to this, the convergence of each method is justified under reasonable assumptions. Gradient based techniques for smooth multi-objective optimization problems are developed in \cite{mat1,mat2,mat4,mat5,flg1,flg2,flg3}. Major contribution in this area is the Newton method for unconstrained multi-objective optimization problems by Fliege et al. (\cite{flg1}). This method uses a quadratic approximation of each objective function to find a descent direction at any iterating point. It is justified that this method converges quadratically under some mild assumptions. The idea of Newton method in \cite{flg1} is extended in \cite{flg3,mat4,mat5} to constrained multi-objective optimization problems. Apart from gradient based techniques, recently some new techniques are developed for nonsmooth multi-objective optimization problems (see\cite{bello1,bento1,bonnel1,montonen1,neto1,qu2,tanabe1}) by several researchers. This techniques are possible extension of the single objective subgradient method (\cite{bello1,montonen1,neto1}), proximal point method (\cite{bento1,bonnel1}), and proximal gradient method (\cite{tanabe1}) etc. to multi-objective case. \\

Apart from convergence, spreading of approximate Pareto front is a major issue for solving multi-objective optimization problems. In \cite{flg1,flg2,mat2,mat4}, multi-start technique is used to generate approximate Pareto front. New initial point selection techniques are developed in \cite{flg3,mat5} to ensure the spreading of Pareto front. Continuation methods for nonlinear multi-objective optimization problems are studied in \cite{sdd1,hillermeier1,martinb1}. In continuation methods, given an initial set of KKT-points ($\mathcal{S}_0$), all further solutions are computed. However, the solutions generated by these methods are restricted to the connected components of the set of KKT-points which contain a point $s\in\mathcal{S}_0$. Other major limitations of these methods are the Hessians requirement and the lack of strategies to handle inequality constraints. To overcome these limitations, predictor corrector methods are developed in \cite{martin1,martin2} based on interval analysis and parallelotope domains. Some other different approaches are designed in \cite{wang2013zigzag,pereyra20131} to generate well distributed approximate Pareto fronts. \\ 

Proximal gradient methods are considered as efficient techniques to solve composite single objective optimization problems (see \cite{beck1,beck2,nesterov1}). At every iteration of this method linear approximation of smooth function is used to find a suitable descent direction. The ideas of proximal gradient methods are further extended by several researchers in various directions. The proximal gradient method developed by Lee et al. (\cite{lee1}) uses a quadratic approximation of smooth function in every iteration. This method converges quadratically under some mild assumptions.\\

Recently Tanabe et al. (\cite{tanabe1}) have developed a proximal gradient method for multi-objective optimizations. This method combines the ideas of the steepest descent method and the proximal point method for multi-objective optimization problems developed in \cite{flg2} and \cite{bonnel1} respectively. At every iteration of this method, a subproblem is solved to find a descent direction which uses linear approximation each smooth function. Similar to the single objective proximal gradient methods, the convergence rate of this method is low. In this paper, we have adopted some ideas of \cite{lee1} and developed a Newton-type proximal gradient method for composite multi-objective optimization problems. The proposed method combines the ideas of Newton method and proximal point method for multi-objective optimization problems developed in \cite{flg1} and \cite{bonnel1} respectively.\\

The outline of the paper is as follows. Some prerequisites are discussed in Section \ref{secpre}. A Newton-type proximal gradient is developed in Section \ref{nprox_method}. An algorithm is proposed in Section \ref{secalg}. The global convergence of this algorithm is justified in this section. In Section \ref{sec_ex}, the proposed method is verified and compared with some existing methods using a set of problems.
\section{Preliminaries}
\lb{secpre}
Consider the multi-objective optimization problem:
\bee
(MOP):~~\underset{x\in\mathbb{R}^n}{\min}~~ F(x)=(F_1(x),F_2(x),...,F_m(x)).\\
\eee
Suppose $F_j:\mathbb{R}^n\rightarrow \mathbb{R}$ is defined by $F_j(x)=f_j(x)+g_j(x)$ where $f_j$ is convex and continuously differentiable and $g_j$ is convex and continuous but not necessarily differentiable function, for $j=1,2,...,m$. 
Denote $\Lambda_n=\{1,2,...,n\}$ for any $n\in \mathbb{N}$. Inequality in $\mathbb{R}^m$ is understood component wise. If there exists $x\in \mathbb{R}^n$ such that $x$ minimizes all objective functions simultaneously, then it is an ideal solution. But in practice, a decrease of one objective function may cause an increase of another objective function. So in the theory of multi-objective optimization optimality is replaced by efficiency. A point $x^*\in \mathbb{R}^n$ is said to be an efficient solution of the $(MOP)$ if there does not exist $x\in \mathbb{R}^n$ such that $F(x)\leq F(x^*)$ and $F(x)\neq F(x^*)$. A feasible point $x^*\in \mathbb{R}^n$ is said to be a weak efficient solution of the $(MOP)$ if there does not exist $x\in \mathbb{R}^n$ such that $F(x)< F(x^*)$. It is clear that every efficient solution of the $(MOP)$ is a weak efficient solution, but the converse is not true. If each $F_j$, $j\in\Lambda_m$ are strictly convex then every weak efficient solution is an efficient solution. For $x, y \in \mathbb{R}^n$, we say $y$ dominates $x$, if and only if $F(y)\leq F(x)$, $ F(y)\neq F(x)$. If $X^*$ is the set of all efficient solutions of the $(MOP)$, then $F(X^*)$ is said to be the Pareto front of the $(MOP)$ and it lies on the boundary of $F(\mathbb{R}^n)$.\\

Define $\mathbb{R}^m_{+}:=\{x\in\mathbb{R}^m\mid x_i\geq 0~\forall~i\in\Lambda_m\}$ and $\mathbb{R}^m_{++}:=int(\mathbb{R}^m_{+})$. Suppose $x^*$ be a weak efficient solution of the $(MOP)$. Then $x^*$ must satisfy $$\left(F_1^{'}(x^*;d),F_2^{'}(x^*;d),...,F_m^{'}(x^*;d)\right)\notin -\mathbb{R}^m_{++}$$ for all $d\in \mathbb{R}^n$. This shows that
\bea
\underset{j\in\Lambda_m}{\max} F_j^{'}(x^*;d)\geq 0 ~~~~\mbox{ for all $d\in\mathbb{R}^n$} \lb{A*}
\eea
The inequality in (\ref{A*}) is sometimes refereed to in the literature as the criticality condition or the first order necessary for weak efficiency of the $(MOP)$ and $x^*$ satisfying (\ref{A*}) is often called a critical point for the $(MOP)$. Further convexity of each $F_j$, $j\in \Lambda_m$ ensures that every critical point of the $(MOP)$ is a weak efficient solution. Further note that if each $F_j$, $j\in\Lambda_m$ is a strictly or strongly convex function, then the critical point of the $(MOP)$ is an efficient solution. Note that if either $f_j$ or $g_j$ is strictly or strongly convex then is so $F_j$. A function $f:\mathbb{R}^n \rightarrow\mathbb{R}$ is said to be a $\sigma$- strongly convex if for every $x,y \in \mathbb{R}^n $ and $\alpha\in [0,1]$,
$$f\left(\alpha x+(1-\alpha)y\right)\leq \alpha f(x)+(1-\alpha)f(y)-\frac{1}{2}\sigma \alpha (1-\alpha)\|x-y\|^2.$$
If $f$ is twice differentiable, then $d^T\nabla^2 f(x) d\geq \sigma \|d\|^2$ for all $d\in\mathbb{R}^n$.\\

In nonsmooth optimization, the concept of the gradient (in smooth optimization) is replaced by the subdifferential. It plays an important role in nonsmooth optimization. The subdifferential of a continuous convex function is defined as follows.
\begin{d1}\lb{subgrad0}(\cite{jd0})
Suppose $h:\mathbb{R}^n\rightarrow (-\infty,\infty]$ be a proper function and $x\in$ $dom~h$. Then
subdifferential of $h$ at $x$ is denoted by $\partial h(x)$ and defined as
\bee
\partial h(x):=\{\xi\in\mathbb{R}^n|h(y)\geq h(x)+\xi^T(y-x)\mbox{ for all }y\in \mathbb{R}^n\}.
\eee
If $x\notin dom~h$ then we define $ \partial h(x)=\emptyset$.
\end{d1}
The following properties of subdifferential are often used in the derivation of the proposed methodology.
\begin{t1}\lb{subgrad1}
({\bf Proposition 2.82, \cite{jd0}}) Let $h:\mathbb{R}^n\rightarrow (-\infty,\infty]$  be a proper convex function, and assume that $x\in int(dom~h)$. Then $ \partial h(x)$ is nonempty and bounded. Moreover, if $h$ is continuous at $x\in dom~h$, then $\partial h (x)$ is compact.
\end{t1}
\begin{t1}\lb{subgrad2}
({\bf Theorem 2.91, \cite{jd0}})
Consider two proper convex functions $h_i:\mathbb{R}^n\rightarrow \mathbb{R}$, $i=1,2$ Suppose that $ri~ dom~h_1 \cap  ri~ dom ~h_2 \neq \emptyset$. Then
\bee
\partial (h_1(x)+h_2(x))=\partial h_1(x)+\partial h_2(x)
\eee
for every $x \in dom~h_1+h_2$.
\end{t1}
\begin{t1}\lb{subgrad3}
({\bf Theorem 2.96, \cite{jd0}}) Consider convex functions $h_j:\mathbb{R}^n\rightarrow \mathbb{R}$, $j= 1,2, \dots,m$ and let $h(x)=\max \{ h_1(x),h_2(x),\dots,h_m(x)\}$. Then
\bee
\partial h(x)=Co \underset{ j\in I(x)}{\cup} \partial h_j(x) 
\eee
where $I(x)=\{j\in\Lambda_m\mid h(x)=h_j(x)\}$ is the active index set and $Co$ is the convex hull.
\end{t1}
It is simple to observe that $x^*={arg~\min}_{x\in\mathbb{R}^n} h(x)$ if and only if $0\in\partial h(x^*).$ 
\section{A Newton-type proximal gradient method for the $(MOP)$}\lb{nprox_method}
In this section a descent algorithm is developed for the $(MOP)$ motivated by Newton-type proximal gradient method for scalar convex optimization problems. Prior to that, we prove the following result which is used in the derivation of the proposed method.
\begin{l1} \lb{ness1}
Suppose $0\in Co \underset{ j\in \Lambda_m}{\cup} \partial F_j(x^*)$ for some $x^*\in\mathbb{R}^n$, then $x^*$ is a critical point of the $(MOP)$. 
\end{l1}
{\bf Proof:} For $x^*\in \mathbb{R}^n$, define $\bar{F}(x)=\max_{j\in\Lambda_m} F_j(x)-F_j(x^*)$. Then from Theorem \ref{subgrad3}, $\partial \bar{F}(x^*)=Co \underset{ j\in \Lambda_m}{\cup} \partial F_j(x^*).$ Thus $0\in Co \underset{ j\in \Lambda_m}{\cup} \partial F_j(x^*)$ implies $x^*={arg~\min}_{x\in\mathbb{R}^n} \bar{F}(x)$. Next we show that, $x^*$ satisfies (\ref{A*}). On contrary, suppose there exists $d\in\mathbb{R}^n$ such that $F_j^{\prime}(x^*,d)<0$ for all $j$. Then we can find some $\alpha>0$ sufficiently small such that $F_j(x^*+\alpha d)<F_j(x^*)$ for all $j$. This implies, $\bar{F}(x^*+\alpha d)<0=\bar{F}(x^*)$ holds for some $x^*+\alpha d\in\mathbb{R}^n$. This contradicts that $x^*={arg~\min}_{x\in\mathbb{R}^n} \bar{F}(x)$. Hence $x^*$ satisfies (\ref{A*}). This implies $x^*$ is a critical point of the $(MOP)$.\qed
$ $\\

Next, we construct a subproblem at $x \in\mathbb{R}^n$ to find a descent direction of the $(MOP)$. Using the ideas of Newton method in \cite{flg1}, we use the following approximation of $F_j$ at $x.$
\bee
Q_{j}(x,d):=\nabla f_j(x)^Td+\frac{1}{2} d^T\nabla^2 f_j(x) d+ g_j(x+d)-g_j(x).
\eee
Define $Q(x,d):=\underset{j\in\Lambda_m}{\max} Q_{j}(x,d).$ For any fixed $x$, $ Q_{j}$ is continuous $d$ and hence $Q$ is continuous $d$. From Theorem \ref{subgrad1}, $\partial_{d} Q_j(x,d)$ is nonempty and bounded. Denote $ I(x,d):=\{j\in\Lambda_m|Q(x,d)= Q_{j}(x,d)\}.$
For any $x\in\mathbb{R}^n$, we solve the following subproblem to find a suitable descent direction of the $(MOP)$.
\bee
P(x):~~\underset{d\in\mathbb{R}^n}{\min} Q(x,d).
\eee 
\begin{n1}
Once can observe that, if $g_j=0$ for all $j\in\Lambda_m$, the $P(x)$ coincides with the subproblem in \cite{flg1} and if $m=1$ then $P(x)$ coincides with the subproblem used in Step 3 of Algorithm 1 in \cite{lee1}.
\end{n1}
If $f_j$ is $\sigma$-strongly convex function for all $j$, then $Q_j(x,d)$ is $\sigma$-strongly convex function in $d$ for every $j\in\Lambda_m$ and $x\in \mathbb{R}^n$. Hence $Q(x,d)$ is $\sigma$-strongly convex function in $d$ for every $x\in \mathbb{R}^n.$ This implies, $P(x)$ has a unique finite minimizer.
\\
Denote $d(x)=\underset{d\in\mathbb{R}^n}{arg~\min}~ Q(x,d)$ and $t(x)=Q(x,d(x)).$ Clearly for every 
$x\in \mathbb{R}^n,$
\bea
t(x)=Q(x,d(x))\leq Q(x,0)=0.\lb{tt0} 
\eea
Since $d(x)$ is the solution of $P(x)$,
\bee
0\in \partial_{d} Q(x,d(x))
\eee
Then from Theorem \ref{subgrad3}, there exists $\lambda\in\mathbb{R}^{|I(x,d(x))|}_{+}$ and $\xi_j\in \partial_{d} g_j(x,d(x))$ $j\in I(x,d(x))$ such that the following conditions hold:
\bee
\underset{j\in I(x,d(x))}{\sum}\lambda_j&=&1\\
\underset{j\in I(x,d(x))}{\sum} \lambda_j\left(\nabla f_j(x)+\nabla^2 f_j(x) d(x)+\xi_j\right)&=&0.
\eee
Substituting $\lambda_j=0$ and using any $\xi_j\in  \partial_{d} g_j(x+d(x))$ for all $j\notin I(x,d(x))$, we can write
\bea
\underset{j\in \Lambda_m}{\sum}\lambda_j&=&1\lb{kkt1}\\
\underset{j\in \Lambda_m}{\sum} \lambda_j\left(\nabla f_j(x)+\nabla^2 f_j(x) d(x)+\xi_j\right)&=&0\lb{kkt2}\\
\lambda_j\geq 0,~\lambda_j\biggl(\nabla f_j(x)^Td(x)+\frac{1}{2} d(x)^T\nabla^2 f_j(x) d(x)\biggr.~~~~~~~\nn\\
\biggl. + g_j(x+d(x))-g_j(x)-t(x)\biggr)&=&0~~j\in\Lambda_m\lb{kkt3}\\ 
\nabla f_j(x)^Td(x)+\frac{1}{2} d(x)^T\nabla^2 f_j(x) d(x)+ g_j(x+d(x))-g_j(x)&\leq& t(x)~~j\in\Lambda_m.\nn\\ \lb{kkt4} 
\eea
Thus if $d(x)$ is the solution of $P(x)$ and $t(x)=Q(x,d(x))$ then there exists $\lambda\in\mathbb{R}^m_{+}$ such that $(d(x),t(x),\lambda)$ satisfies (\ref{kkt1})-(\ref{kkt4}).
\begin{l1}\lb{lm1}
Suppose $f_j$ is strictly convex function for all $j$. Then $x\in\mathbb{R}^n$ is a critical point of the $(MOP)$ if and only if $d(x)=0.$ 
\end{l1}
{\bf Proof:} If possible let $x$ is a critical point of the $(MOP)$ and $d(x)\neq0$. Since $f_j$ is strictly convex for every $j$, from (\ref{tt0})
 \bea
 \nabla f_j(x)^Td(x)+g_j(x+d(x))-g_j(x)\leq -\frac{1}{2}d(x)^T \nabla^2 f_j(x) d(x)<0.\lb{lm100}
 \eea
Since $g_j$ is convex, for any $\alpha\in (0,1)$ 
\bea
g_j(x+\alpha d(x))-g_j(x)&\leq&\alpha g_j(x+d(x))+(1-\alpha)g_j(x)-g_j(x)\nn\\
&=&\alpha\left(g_j(x+d(x))-g_j(x)\right). \lb{lm111}
\eea
From (\ref{lm100}) and (\ref{lm111}),
\bee
& &\alpha \nabla f_j(x)^Td(x)+g_j(x+\alpha d(x))-g_j(x)\\
&\leq& \alpha \left(\nabla f_j(x)^Td(x)+g_j(x+d(x))-g_j(x)\right)\\
&<&0.
\eee 
This implies
\bee
\frac{1}{\alpha}\left(\alpha \nabla f_j(x)^Td(x)+g_j(x+\alpha d(x))-g_j(x)\right)<0
\eee
Taking limit $\alpha\rightarrow 0^{+}$ in the above inequality we have $F_{j}^{\prime}(x,d(x))<0$ for all $j\in \Lambda_m.$ This shows that $x$ is not a critical point, a contradiction. Hence, if $x$ is a critical point, then $d(x)=0.$\\
Conversely suppose $d(x)=0$. Then from (\ref{kkt1}) and (\ref{kkt2}), there exists $\lambda\in\mathbb{R}^m_{+}$ such that $\sum_{j\in\Lambda_m} \lambda_j=1$ and $\sum_{j\in \Lambda_m} \lambda_j\left(\nabla f_j(x)+\xi_j\right)=0$ where $\xi_j\in\partial g_j(x)$ for $j\in\Lambda_m$. This implies $$0\in Co \underset{ j\in \Lambda_m}{\cup} \partial F_j(x).$$       
Hence, from Lemma \ref{ness1}, $x$ is a critical point of $(MOP).$\qed
\begin{n1}
If $f_j$ is strictly convex for all $j$ then $P(x)$ has a unique solution. This implies $t(x)=0$ holds if and only if $d(x)=0$. Hence, from Lemma \ref{lm1} and (\ref{tt0}) we can conclude that $t(x)<0$ holds if and only if $x$ is a non-critical point of the $(MOP)$.
\end{n1}
\begin{t1}\lb{th11}
Suppose, for $x^k\in\mathbb{R}^n$ $d(x^k)$ be the solution of $P(x^k)$. Further suppose, $\{x^k\}$ converges to $x^*$ as $k\rightarrow \infty$ and the assumptions of Lemma \ref{lm1} hold. If $\{d(x^k)\}$ converges to $d^{*}$ then $d^{*}=d(x^{*}).$  
\end{t1}
{\bf Proof:} Suppose, for $x^k\in \mathbb{R}^n$, $d(x^k)$ be the solution of $P(x^k)$ and $t(x^k)=Q(x^k,d(x^k)).$ Then from (\ref{kkt1})-(\ref{kkt4}), there exists $\lambda^k\in\mathbb{R}^m_{+}$ such that 
\bea
\underset{j\in \Lambda_m}{\sum}\lambda_j^k&=&1\lb{kkkt1}\\
\underset{j\in \Lambda_m}{\sum} \lambda_j^k\left(\nabla f_j(x^k)+\nabla^2 f_j(x^k) d(x^k)+\xi_j^k\right)&=&0\lb{kkkt2}\\
\lambda_j^k\geq 0,~\lambda_j^k\biggl(\nabla f_j(x^k)^Td(x^k)+\frac{1}{2} d(x^k)^T\nabla^2 f_j(x^k) d(x^k)\biggr.& &\nn\\
\biggl.+ g_j(x^k+d(x^k))-g_j(x^k)-t(x^k)\biggr)&=&0~~j\in\Lambda_m\lb{kkkt3}\\ 
\nabla f_j(x^k)^Td(x^k)+\frac{1}{2} d(x^k)^T\nabla^2 f_j(x^k) d(x^k)& &\nn\\+ g_j(x^k+d(x^k))-g_j(x^k)&\leq& t(x^k)~~j\in\Lambda_m \lb{kkkt4}.
\eea 
Since $$\partial_{d} \left(g_j(x^k+.)-g_j(x^k)\right)(d(x^k))=\partial_{d} g_j(x^k+d(x^k))$$ is bounded for every $j$, the sequence $\{(\xi_1^k,\xi_2^k,...,\xi_m^k)\}$ is bounded. Then there exists a convergent subsequence $\{(\xi_1^k,\xi_2^k,...,\xi_m^k)\}_{k\in K_1}$ converging to $(\xi_1^*,\xi_2^*,...,\xi_m^*).$ From (\ref{kkkt1}), $\{\lambda^k\}_{k\in K_1}$ is bounded and hence there exists a sub subsequence $\{\lambda^k\}_{k\in K_2\subseteq K_1}$ converging to $\lambda^*.$ Now taking sum over $j\in\Lambda_m$ in (\ref{kkkt3}) and using (\ref{kkkt1}),
\bee
\underset{j\in \Lambda_m}{\sum}\lambda_j^k\biggl(\nabla f_j(x^k)^Td(x^k)+\frac{1}{2} d(x^k)^T\nabla^2 f_j(x^k) d(x^k)\biggr.\\\biggl.+ g_j(x^k+d(x^k))-g_j(x^k)\biggr)=t(x^k)
\eee
Suppose
\bea
t^{*}:&=&\underset{\underset{k\in K_2}{k\rightarrow \infty}}{\lim} t(x^k) \nn\\
&=& \underset{j\in \Lambda_m}{\sum}\lambda_j^*\biggl(\nabla f_j(x^*)^Td(x^*)+\frac{1}{2} d(x^*)^T\nabla^2 f_j(x^*) d(x^*)\biggr.\nn\\
\biggl. & & +g_j(x^*+d(x^*))-g_j(x^*)\biggr) \lb{*1}.
\eea
Considering limit $k\rightarrow \infty$, $k\in K_2$ in (\ref{kkkt1})-(\ref{kkkt4}),
\bea
\underset{j\in \Lambda_m}{\sum}\lambda_j^*&=&1\lb{skkt1}\\
\underset{j\in \Lambda_m}{\sum} \lambda_j^*\left(\nabla f_j(x^*)+\nabla^2 f_j(x^*) d^*+\xi_j^*\right)&=&0\lb{skkt2}\\
\lambda_j^*\geq 0,~\lambda_j^*\biggl(\nabla f_j(x^*)^Td^*+\frac{1}{2} {d^*}^T\nabla^2 f_j(x^*) d^*\biggr.& &\nn\\ \biggl.+ g_j(x^*+d^*)-g_j(x^*)-t^*\biggr)&=&0~~j\in\Lambda_m \lb{skkt3}\\ 
\nabla f_j(x^*)^Td^*+\frac{1}{2} {d^*}^T\nabla^2 f_j(x^*) d^*+ g_j(x^*+d^*)-g_j(x^*)&\leq& t^*~~j\in\Lambda_m\nn\\ \lb{skkt4}.
\eea 
Since right hand side of (\ref{skkt4}) is independent of $j$, 
\bea
\underset{j\in\Lambda_m}{\max} \nabla f_j(x^*)^Td^*+\frac{1}{2} {d^*}^T\nabla^2 f_j(x^*) d^*+ g_j(x^*+d^*)-g_j(x^*)\leq t^*.\lb{lm21}
\eea
From (\ref{*1}) and (\ref{skkt1}),
\bea
\underset{j\in\Lambda_m}{\max} \nabla f_j(x^*)^Td^*+\frac{1}{2} {d^*}^T\nabla^2 f_j(x^*) d^*+ g_j(x^*+d^*)-g_j(x^*)\geq t^*.\lb{lm22}
\eea
Hence from (\ref{lm21}) and (\ref{lm22}),
\bee
\underset{j\in\Lambda_m}{\max} \nabla f_j(x^*)^Td^*+\frac{1}{2} {d^*}^T\nabla^2 f_j(x^*) d^*+ g_j(x^*+d^*)-g_j(x^*)= t^*
\eee
Denote $$\tilde{I}(x^*,d^*)=\{j\in\Lambda_m|\nabla f_j(x^*)^Td^*+\frac{1}{2} {d^*}^T\nabla^2 f_j(x^*) d^*+ g_j(x^*+d^*)-g_j(x^*)= t^*\}.$$
Then from (\ref{skkt3}), $\lambda_j^*=0$ for all $j\notin \tilde{I}(x^*,d^*).$ Since $\xi_j^k\in\partial_{d}g_j(x^k+d^k)$,
\bee
g_j(x^k+d)-g_j(x^k+d(x^k))\geq {\xi_j^k}^T(d-d(x^k))
\eee
holds for all $d\in \mathbb{R}^n.$ Taking limit $k\rightarrow \infty$, $k\in K_2$, 
\bee
g_j(x^*+d)-g_j(x^*+d^*)\geq {\xi_j^*}^T(d-d^*)
\eee
for all $d\in \mathbb{R}^n$. This implies $\xi_j^*\in \partial_{d}g_j(x^*+d^*)$. Hence from (\ref{skkt1}) and (\ref{skkt2}),
\bee
0\in Co\underset{j\in\Lambda_m}{\cup} \left\{\nabla f_j(x^*)+\nabla^2f_j(x^*)d^*\right\}+ \partial_{d}g_j(x^*+d^*)
\eee 
This implies $d^*$ is a solution of $P(x^*)$. Since $Q(x^*,d)$ is strictly convex in $d$, $P(x^*)$ has unique solution. Hence $d^*=d(x^*)$ and $t^*=t(x^*)$.\qed  
$ $\\

Next we develop an Armijo type line search technique to find a suitable step length that ensures sufficient decrease in each objective function. We consider $\alpha>0$ as a suitable step length if for some $\beta\in (0,1)$
\bea
F_j(x+\alpha d(x))\leq F_j(x)+\beta\alpha t(x) \lb{amj1}
\eea 
holds for every $j\in\Lambda_m.$\\ 

In the following theorem we have justified that the above line search technique is well defined.
\begin{t1}\lb{t1}
Suppose $f_j$ is $\sigma$-strongly convex for all $j$ and $d(x)$ is the solution of $P(x)$. Then 
\bea
t(x)\leq -\frac{\sigma}{2}\|d(x)\|^2.\lb{t10}
\eea 
Further, if $x$ is non critical then (\ref{amj1}) holds for every $\alpha>0$ sufficiently small.
\end{t1}
\textbf{Proof:} Suppose $d(x)$ is the solution of $P(x)$ and $t(x)=Q(x,d(x))$. Then there exists $\lambda\in\mathbb{R}^m_{+}$ such that $(d(x),t(x),\lambda)$ satisfies (\ref{kkt1})-(\ref{kkt4}). Since $g_j$ is convex and $\xi_j\in\partial_{d}g_j(x+d(x)),$
\bea
g_j(x+d(x))-g_j(x)\leq \xi_j^T d(x).\lb{t11}
\eea
Multiplying both sides of (\ref{kkt2}) by $d(x)$,
\bee
\underset{j\in \Lambda_m}{\sum} \lambda_j\left(\nabla f_j(x)^Td(x)+d(x)^T\nabla^2 f_j(x) d(x)+\xi_j^Td(x)\right)=0
\eee
Hence from (\ref{t11}),
\bea
\underset{j\in \Lambda_m}{\sum} \lambda_j\left(\nabla f_j(x)^Td(x)+d(x)^T\nabla^2 f_j(x) d(x)+g_j(x+d(x))-g_j(x)\right)\leq 0.\lb{t12}
\eea
Taking sum over $j\in\Lambda_m$ in (\ref{kkt3}) and using (\ref{kkt1}),
\bea
& &\underset{j\in\Lambda_m}{\sum}\lambda_j\left(\nabla f_j(x)^Td(x)+ d(x)^T\nabla^2 f_j(x) d(x)+ g_j(x+d(x))-g_j(x)\right)\nn\\
&= &\underset{j\in\Lambda_m}{\sum}\lambda_j\frac{1}{2}d(x)^T\nabla^2 f_j(x) d(x)+t(x).\lb{t13}
\eea 
Using (\ref{t12}) in (\ref{t13}),
\bea
t(x)\leq - \underset{j\in\Lambda_m}{\sum}\lambda_j\frac{1}{2}d(x)^T\nabla^2 f_j(x) d(x).\lb{t14}
\eea 
Since $f_j$ is $\sigma$-strongly convex for every $j$, $d(x)^T\nabla^2 f_j(x) d(x)\geq \sigma\|d(x)\|^2$ holds for every $j$. Then from (\ref{t14}) and (\ref{kkt1}),
\bee
t(x)\leq-\frac{\sigma}{2}\|d(x)\|^2. 
\eee
Suppose $x$ is non critical. Then from Lemma \ref{lm1}, $d(x)\neq 0.$ Hence, from (\ref{t10}), $t(x)<0.$ Since $g_j$ is convex, for any $\alpha\in[0,1]$,
\bee
F_j(x+\alpha d(x))-F_j(x)&=&f_j(x+\alpha d(x))-f_j(x)+g_j(x+\alpha d(x))-g_j(x)\\
&\leq&\alpha\left(\nabla f_j(x)^Td(x)+g_j(x+d)-g_j(x)\right)+\mathcal{O}(\alpha^2)\\
&<&\alpha t(x) +\mathcal{O}(\alpha^2)~~~~(\because~~\frac{1}{2}d(x)^T\nabla^2f_j(x)d(x)>0)
\eee 
Then for every $j\in\Lambda_m$, 
\bea
F_j(x+\alpha d(x))-F_j(x)-\alpha \beta t(x) <\alpha(1-\beta)t(x)+\mathcal{O}(\alpha^2)\lb{t15}
\eea
Since $\beta\in(0,1)$ and $t(x)<0$, the right hand side term in  (\ref{t15}) becomes non positive for every $\alpha>0$ sufficiently small. This implies, (\ref{amj1}) holds for every $\alpha>0$ sufficiently small.\\
Hence the theorem follows.\qed
\section{Algorithm and convergence analysis}\lb{secalg}
In this section we develop an algorithm for the $(MOP)$ using the theoretical results developed so far. In addition to this, we justify the global convergence of this algorithm under some mild assumptions. For simplicity, the rest of the paper, we denote $d(x^k)$ and $t(x^k)$ by $d^k$ and $t^k$ respectively.
\begin{alg1} \lb{npxm1}(Newton-type proximal gradient method for the $(MOP)$)
\begin{enumerate}[{Step} 1]
\item Choose initial approximation $x^0$, scalars $r,\beta \in(0,1)$, and $\epsilon>0$.Set $k:=0$
\item Solve the subproblem $P(x^k)$ to find $d^k$ and $t^k$.\lb{P_S}
\item If $\|d^k\|<\epsilon$, then stop. Else go to Step \ref{st_alp}. \lb{st_term}
\item Choose a suitable step length $\alpha_k$ as the first element of $\{1,r,r^2,...\}$ satisfying (\ref{amj1}).\lb{st_alp}
\item Update $x^{k+1}=x^k+\alpha_k d^k$, set $k:=k+1$, and go to Step \ref{P_S}.
\end{enumerate}
\end{alg1}
Global convergence of this algorithm is justified in the following theorem. 
\begin{t1}\lb{cv0}
Suppose $\{x^k\}$ is a sequence generated by Algorithm \ref{npxm1}, $f_j$ is $\sigma$-strongly convex functions for every $j$, $\nabla f_j$ is Lipschitz continuous for every $j$ with Lipschitz constant $L$, and the level set $M:=\{x\in\mathbb{R}^n|F(x)\leq F(x^0)\}$ is bounded. Then every accumulation point of  $\{x^k\}$ is a critical point of the $(MOP)$.
\end{t1}
{\bf Proof:} First, we show that there exists $\bar{\alpha}>0$ such that $\alpha_k\geq \bar{\alpha}$ holds for every $k$. Since $\nabla f_j$ is Lipschitz continuous, for any $\alpha$ we have
\bea
f_j(x^k+\alpha d^k)\leq f_j(x^k)+\alpha \nabla f_j(x^k)^Td^k+\frac{L}{2}\alpha^2\|d^k\|^2 ~~
\forall~~j\in\Lambda_m. \lb{cv01}
\eea
From (\ref{cv01}) and (\ref{lm111}) for any $\alpha\in[0,1]$,
\bea
F_j(x^k+\alpha d^k)-F_j(x^k)&\leq &\alpha \left(\nabla f_j(x^k)^Td^k+g_j(x^k+d^k)-g_j(x^k)\right)+\frac{L}{2}\alpha^2\|d^k\|^2\nn\\
&<& \alpha \left(\nabla f_j(x^k)^Td^k+\frac{1}{2} {d^k}^T\nabla^2 f_j(x^k) d^k+g_j(x^k+d^k)-g_j(x^k)\right)\nn\\
& &~~+\frac{L}{2}\alpha^2\|d^k\|^2\nn\\ 
&\leq &\alpha t^k+\frac{L}{2}\alpha^2\|d^k\|^2~~\forall~~j\in\Lambda_m \lb{cv03}.
\eea
Second inequality holds since ${d^k}^T\nabla^2 f_j(x^k) d^k>0$. From Step \ref{st_alp} of Algorithm \ref{npxm1}, either $\alpha_k=1$ or there exists $k_1\in \mathbb{N}$ such that $\alpha_k=r^{k_1}$. If $\alpha_k=r^{k_1}$ then there exists at least one $\hat{j}\in\Lambda_m$ such that
\bee
F_{\hat{j}}(x^k+r^{k_1-1}d^k)-F_{\hat{j}}(x^k)>r^{k_1-1}\beta t^k.
\eee
Then from (\ref{cv03}), $r^{k_1-1}\beta t^k<r^{k_1-1} t^k+\frac{L}{2}r^{2(k_1-1)}\|d^k\|^2.$
This implies
\bea
-(1-\beta)t^k&<&\frac{L}{2}r^{(k_1-1)}\|d^k\|^2 \lb{cv04}.
\eea
Using (\ref{t10}) in (\ref{cv04}), $\frac{\sigma(1-\beta)}{2}\|d^k\|^2<\frac{L}{2}r^{k_1-1}\|d^k\|^2.$ This implies $r^{k_1}>\frac{\sigma(1-\beta)r}{ L}.$ Choose $\bar{\alpha}=\min\{1,\frac{\sigma(1-\beta)r}{L}\}$. Then $\alpha_k\geq\bar{\alpha}$ holds for every $k$.\\

Now from Step \ref{st_alp} of Algorithm \ref{npxm1} for any $N\in\mathbb{N}$ and $j\in\Lambda_m$
\bea
F_j(x^{N+1})-F_j(x^0)&\leq&\beta\sum_{k=0}^{N} \alpha_kt^k\nn\\
&\leq&-\frac{\sigma\beta}{2}\sum_{k=0}^{N} \alpha_k \|d^k\|^2\nn\\
&\leq&-\frac{\sigma\beta\bar{\alpha}}{2}\sum_{k=0}^{N}  \|d^k\|^2 \lb{cv06}
\eea
Second inequality follows from (\ref{t10}) and last inequality holds since $\alpha_k\geq \bar{\alpha}.$ From (\ref{amj1}) and (\ref{t10}), $x^k\in M$ for every $k$. That is $\{x^k\}$ is bounded. Taking limit $N\rightarrow \infty$ in (\ref{cv06}) and using the continuity of $F_j$,
\bee
-\frac{\sigma\beta\bar{\alpha}}{2}\sum_{k=0}^{\infty}  \|d^k\|^2\geq \underset{k\rightarrow\infty}{\lim} F_j(x^k)-F_j(x^0)>-\infty.
\eee
This implies $\underset{k\rightarrow \infty}{\lim} d^k=0.$ Since $\{x^k\}$ is bounded, there exists convergent subsequence $\{x^k\}_{k\in K}$ converging to $x^*$. Since $\underset{k\rightarrow \infty}{\lim} d^k=0$, $\underset{\underset{k\in K}{k\rightarrow \infty}}{\lim} d^k=0$. Then from Theorem \ref{th11}, $d(x^*)=0$. From Lemma \ref{lm1}, $x^*$ is a critical point of $(MOP).$\\
Hence the theorem follows. \qed
\section{Numerical examples}
\lb{sec_ex}
In this section Algorithm \ref{npxm1} (MONPG) is verified and compared with the proximal gradient method for (MOP) developed in \cite{tanabe1} (MOPG) and weighted sum method (WS) using a set of problems.  Implementation of these algorithms is explained below.
\begin{itemize}
\item MATLAB (2019a) code is developed for each method using the extension `{\it CVX}'. `{\it CVX}' programming with solver `{\it SeDuMi}' is used to solve the subproblems in each method.
\item The single objective optimization problem of the weighted sum method is solved using the Newton-type proximal gradient method developed in \cite{lee1}.
\item $\|d^k\|<10^{-5}$ or maximum 200 iterations is considered as stopping criteria.
\item Solution of a multi-objective optimization problem is not isolated optimum points, but a set of efficient solutions. To generate an approximate set of efficient solutions we have considered multi-start technique. Following steps are executed in this technique.
\begin{itemize}
\item A set of 100 uniformly distributed random initial points between $lb$ and $ub$ are considered, where $lb,ub\in\mathbb{R}^n$ and $lb<ub$.
\item Algorithm \ref{npxm1} is executed individually.
\item Similar process is followed in case of MOPG.
\item For weighted sum method 100 different weights are used. For 2 objective optimization problems, we have used $(1,0)$, $(0,1)$ and 98 uniformly distributed weights. For 3 objective optimization problems $(1,0,0)$, $(0,1,0)$, $(0,0,1)$ and 97 uniformly distributed weights are used.
\item A uniformly distributed random initial approximation is used to solve each weighted single objective optimization problem.
\item Suppose $\mathcal{WX^*}$ is the collection of approximate critical points. The non dominated set of $\mathcal{WX^*}$ is considered as an approximate set of efficient solutions.
\end{itemize} 
\end{itemize}
Next we explain the steps of Algorithm \ref{npxm1} using the following example. \\

\textbf{Example 1:} Consider the bi-objective optimization problem: 
$$(P_1):~~\underset{x\in\mathbb{R}^2}{\min}~~\left(f_1(x)+g_1(x), f_2(x)+g_2(x)\right)$$
where
\bea
f_1(x)&:=& x_1^4+x_2^4,\lb{ex1f1}\\
f_2(x)&:=&(x_1-5)^4+(x_2-5)^4,\lb{ex1f2}\\
g_1(x)&:=& \max\{(x_1-2)^2+(x_2+2)^2,x_1^2+8x_2\},\lb{ex1g1}\\
g_2(x)&:=& \max\{5x_1+x_2,x_1^2+x_2^2\} \lb{ex1g2}.
\eea
(Here $(g_1(x),g_2(x))$ follows from Example 2 in Section 4.1 of \cite{montonen1}.)\\\\
Consider $x^0=(3.7990,1.8743 )^T.$ Then $F(x^0)=(250.0622,118.4027)^T$. Solution of $P(x^0)$ is obtained as $d^0=(-0.6444,0.9601)^T$ and $t^0=-57.4460$. $\alpha_0=1$ satisfies (\ref{amj1}). The next iterating point is $x^1=x^0+\alpha_0d^0=(3.1546,2.8344)^T.$ Clearly $F(x^1)=(196.2014,52.1993)^T<F(x^0)$. Using the stopping condition $\|d^k\|<10^{-5}$, approximate solution is obtained as $x^2=(2.9912,3.0017)^T\approx(3,3)^T.$\\
One can verify that $x^*=(3,3)^T$ is a critical point of $(P_1)$ using weighted sum method with weight $w=(0.18637886,0.81362114).$\\

\begin{figure}[!htbp]
    \centering
    \begin{subfigure}[b]{.47\textwidth}
        \centering
        \includegraphics[height=2.75cm,width=\linewidth]{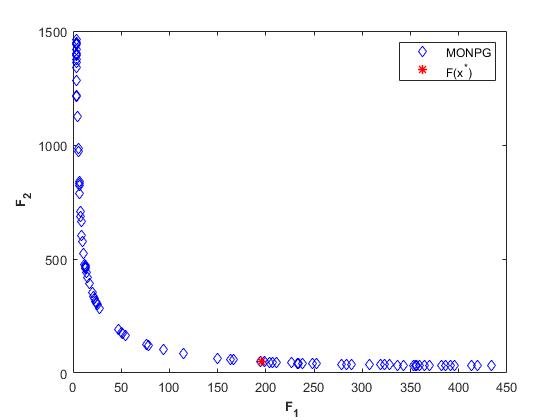}
        \caption{Approximate Pareto fronts by MONPG}
        \label{fig1b}
    \end{subfigure}
    \begin{subfigure}[b]{.47\textwidth}
        \centering
        \includegraphics[height=2.75cm,width=\linewidth]{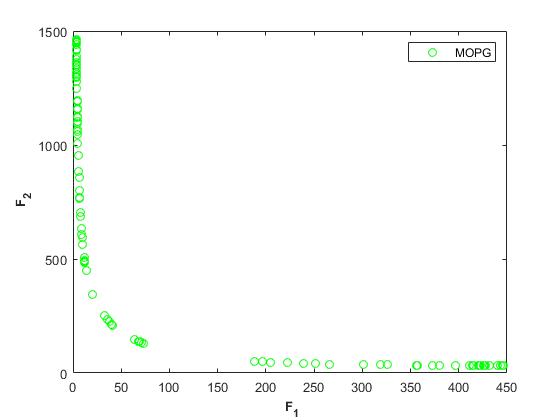}
        \caption{Approximate Pareto fronts by MOPG}
        \label{fig1a}
    \end{subfigure}
\caption{Approximate Pareto fronts of $(P_1)$}
\end{figure}

\begin{figure}[htb]
    \centering
       \includegraphics[height=4cm,width=.8\linewidth]{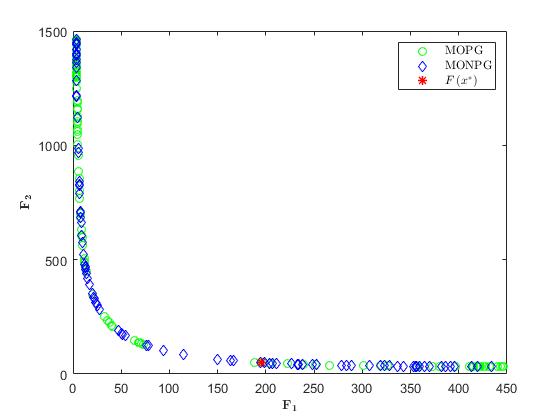}
        \caption{Approximate Pareto fronts of $P_1$ by MONPG and MOPG}
        \label{fig1c}
\end{figure}

Approximate Pareto fronts obtained by multi-start technique with MONPG and MOPG are provided in Figure \ref{fig1a} and \ref{fig1b}. Approximate Pareto fronts are combined in Figure \ref{fig1c}. One can observe that $F(x^*)$ belongs to this approximate Pareto front.
$ $\\\\

{\bf Test Problems:} We have constructed a set of nonlinear 2 or 3 objective optimization problems. Some differentiable multi-objective test problems ($f(x)$) is combined with the following nondifferentibale multi-objective problems. Details of nonsmooth functions are provided here.
\begin{enumerate}[(i).]
\item $gA:\mathbb{R}^2\rightarrow\mathbb{R}^2$ defined in (\ref{ex1g1}) and (\ref{ex1g2}).
\item $gB:\mathbb{R}^2\rightarrow\mathbb{R}^2$ is defined by
\bee
gB_1(x)&=&\max\{x_1^2+(x_2-1)^2,x_1+1)^2\}\\
g_2(x)&=&\max\{x_1^4+x_2^2,2x_1+2x_2\}.
\eee
$gB$ is same as Example 1 in Section 4.1 of \cite{montonen1}.
\item $gC:\mathbb{R}^3\rightarrow\mathbb{R}^3$ is defined by 
\bee
gC_1(x)&=&\max(x_1^2+x_2^2+x_3^2-1,x_1^2+x_2^2+(x_3-2)^2)\\
gC_2(x)&=&\max(x_1+x_2+x_3-1,x_1+x_2-x_3+1)\\
gC_3(x)&=&\max(2x_1^2+6x_2^2+2(5x_3-x_1)^2,x_1^2-9x_3)
\eee
$gC(x)$ is constructed using the ideas of test problem EVD52 in \cite{lukvsan1}. In order to maintain convexity, we have made slight change in the first function of $gC_3(x)$. 
\item $gD:\mathbb{R}^2\rightarrow\mathbb{R}^3$ is defined by 
\bee
gD_1(x)&=&\max(x_1^2+x_2^4,(2-x_1)^2+(2-x_2)^2,2e^{-x_1+x_2})\\
gD_2(x)&=&\max(x_1^4+x_2^2,(2-x_1)^2+(2-x_2)^2,2e^{-x_1+x_2})\\
gD_3(x)&=&\max(5x_1+x_2,-5x_1+x_2,x_1^2+x_2^2+4x_2))
\eee
$gD(x)$ is constructed using the ideas of test problems CB2, CB3, and DEM in \cite{lukvsan1} respectively. 
\item $gE:\mathbb{R}^2\rightarrow\mathbb{R}^3$ is defined by 
\bee
gE_1(x)&=&\max\{x_1^2+x_2^2,x_1^2+x_2^2+10(-4x_1-x_2+4),x_1^2+x_2^2+10(-x_1-2x_2+6)\}\\
gE_2(x)&=&\max\{x_1^2+x_2^4,(2-x_1)^2+(2-x_2)^2,2e^{-x_1+x_2}\}\\
gE_3(x)&=&\max\{5x_1+x_2,-5x_1+x_2,x_1^2+x_2^2+4x_2\}.
\eee
$gE(x)$ is constructed using the ideas of the test problems QL, CB3, DEM in \cite{lukvsan1}.
\item $gF:\mathbb{R}^4\rightarrow\mathbb{R}^2$ is defined by 
\bee
gF_1(x)&=&\max\{g_{11}(x),g_{12}(x)\}\\
gF_2(x)&=&\max\{g_{21}(x),g_{22}(x)\}\\
g_{11}(x)&=&x_1^2+x_2^2+2x_3^2+x_4^2-5x_1-5x_2-21x_3+7x_4\\
g_{12}(x)&=&g_{11}(x)+10(x_1^2+x_2^2+x_3^2+x_4^2+x_1-x_2+x_3-x_4-8)\\
g_{21}(x)&=&g_{11}(x)+10(x_1^2+2x_2^2+x_3^2+2x_4^2-x_1-x_4-10)\\
g_{22}(x)&=&g_{11}(x)+10(2x_1^2+x_2^2+x_3^2+2x_1-x_2-x_4-5).
\eee
$gF(x)$ is constructed using the ideas of the test problem Rosen-Suzuki of \cite{lukvsan1}.
\item $gG:\mathbb{R}^n\rightarrow\mathbb{R}^m$ is defined by 
\bee
g_j(x)&=&\max\{u_{j1}^Tx,u_{j2}^Tx\}~~~j=1,2,\dots,m
\eee
where $u_{j1}$ and $u_{j2}$ $j=1,2,\dots,m$ are uniformly distributed random vectors between  $(0,0,\dots, 0)\in\mathbb{R}^{n}$ and $(0.1,0.1,\dots, 0.1)\in\mathbb{R}^{n}.$
\item $gH:\mathbb{R}\rightarrow\mathbb{R}^2$ is defined by $gH_1(x)=|x|=gH_2(x).$
\end{enumerate}
Details of test problems are provied in table \ref{table1}. 
\begin{table}[!h]
\begin{center}
\small
\begin{tabular}
{|c|c|c|c|c|c|}\hline
Sl. No& $(m, n)$& $f$ & $g$ &$lb^T$ &$ub^T$\\ \hline
\rownumber&(2,2)&(\ref{ex1f1})\&(\ref{ex1f2})& $gA$&$(-3,-3)$& $(7,7)$ \\ \hline
\rownumber&(2,2)&(\ref{ex1f1})\&(\ref{ex1f2})& $gB$&$(-3,-3)$& $(7,7)$\\ \hline
\rownumber&(2,2)&AP2 (\cite{mat1})& gA&$(-5,-5)$&$(5,5)$\\ \hline
\rownumber &(2,2)&AP2 (\cite{mat1})& gB&$(-5,-5)$&$(5,5)$\\ \hline
\rownumber&(2,2)&BK1 (\cite{hub1})& gA&$(-3,-3)$&$(5,5)$\\ \hline
\rownumber &(2,2)&BK1 (\cite{hub1})& gB&$(-3,-3)$&$(5,5)$\\ \hline
\rownumber &(3,3)&FDS (\cite{flg1})& gC&$(-2,-2,-2)$&$(4,4,4)$\\ \hline
\rownumber &(3,5)&FDS (\cite{flg1})& gG&$(-2,\dots,-2)$&$(2,\dots,2)$\\ \hline
\rownumber &(3,8)&FDS (\cite{flg1})& gG&$(-2,\dots,-2)$&$(2,\dots,2)$\\ \hline
\rownumber &(3,2)&IKK1 (\cite{hub1})& gD&$(-2,-2)$&$(3,3)$\\ \hline
\rownumber &(2,2)&Jin1 (\cite{jin1})& gA&$(-3,-3)$&$(5,5)$\\ \hline
\rownumber &(2,2)&Jin1 (\cite{jin1})& gB &$(-3,-3)$&$(5,5)$\\ \hline
\rownumber &(2,4)&Jin1 (\cite{jin1})& gF&$(-5,\dots,-5)^T$&$(10,\dots,10)^T$\\ \hline
\rownumber &(2,4)&Jin1 (\cite{jin1})& gF &$(-100,\dots,-100)^T$&$(100,\dots,100)^T$\\ \hline
\rownumber &(2,10)&Jin1 (\cite{jin1})& gG &$(-5,\dots,-5)^T$&$(5,\dots,5)^T$\\ \hline
\rownumber &(2,10)&Jin1 (\cite{jin1})& gG &$(-100,\dots,-100)^T$&$(100,\dots,100)^T$\\ \hline
\rownumber  &(2,2)&Lovison1 (\cite{lovison1})& gA &$(-3,-3)$&$(5,5)$\\ \hline
\rownumber &(2,2)&Lovison1 (\cite{lovison1})& gB&$(-3,-3)$&$(5,5)$\\ \hline
\rownumber &(2,2)&LRS1 (\cite{hub1})& gB&$(-50,-50)$&$(50,50)$\\ \hline
\rownumber &(2,2)&LRS1 (\cite{hub1})& gA&$(-50,-50)$&$(50,50)$\\ \hline
\rownumber &(3,2)&MHHM1 (\cite{hub1})& gE&$(-4,-4)$&$(4,4)$\\ \hline
\rownumber &(3,2)&MHHM1 (\cite{hub1})& gD &$(-4,-4)$&$(4,4)$ \\ \hline 
\rownumber &(2,1)&MOP1 (\cite{hub1})& gG &$-100$&$100$\\ \hline 
\rownumber &(3,2)&MOP7 (\cite{hub1})& gD &$(-4,-4)$&$(4,4)$\\ \hline 
\rownumber &(3,2)&MOP7 (\cite{hub1})& gE&$(-4,-4)$&$(4,4)$\\ \hline
\rownumber &(2,2)&MS1(\cite{martin1})& gA&$(-2,-2)$&$(2,2)$\\ \hline
\rownumber &(2,2)&MS1 (\cite{martin1})& gB&$(-2,-2)$&$(2,2)$\\ \hline
\rownumber &(2,4)&MS2(\cite{martin1})& gG&$(-2,\dots,-2)$&$(2,\dots,2)$\\ \hline
\rownumber &(2,10)&MS2 (\cite{martin1})& gF&$(-2,\dots,-2)$&$(2,\dots,2)$\\ \hline
\rownumber &(3,3)&SDD1 (\cite{sdd1})& gC&$(-2,-2,-2)$&$(2,2,2)$\\ \hline
\rownumber &(3,10)&SDD1 (\cite{sdd1}) & gC&$(-3,\dots,-3)$&$(3,\dots,3)$\\ \hline
\rownumber&(2,2)&SP1 (\cite{hub1})& gA&$(-1,-1)$&$(5,5)$\\ \hline
\rownumber &(2,2)&SP1 (\cite{hub1})& gB&$(-1,-1)$&$(5,5)$\\ \hline 
\rownumber &(3,2)&VFM1 (\cite{hub1}) & gD&$(-2,-2)$&$(4,4)$\\ \hline
\rownumber &(3,2)&VFM1 (\cite{hub1})& gE &$(-2,-2)$&$(4,4)$ \\ \hline
\rownumber&(2,2)&VU1 (\cite{hub1})& gA&$(-3,-3)$&$(3,3)$\\ \hline
\rownumber &(2,2)&VU1 (\cite{hub1})& gB&$(-3,-3)$&$(3,3)$\\ \hline
\rownumber&(2,2)&VU2 (\cite{hub1})& gA&$(-3,-3)$&$(3,3)$\\ \hline
\rownumber &(2,2)&VU2(\cite{hub1})& gB&$(-3,-3)$&$(3,3)$\\ \hline 
\rownumber &(3,3)&ZLT1 (\cite{hub1})& gC &$(-10,-10,-10)$&$(10,10,10)$\\ \hline 
\rownumber &(3,3)&ZLT1 (\cite{hub1})&gC &$(-100,-100,-100)$&$(100,100,100)$\\ \hline 
\rownumber &(3,10)&ZLT1 (\cite{hub1})& gG &$(-5,\dots,-5)$&$(5,\dots,5)$\\ \hline
\end{tabular}
\caption{Details of test problems}
\lb{table1}
\end{center}
\end{table}  
\\

\textbf{Comparison with existing methods:} MONPG is compared with MOPG and WS. Performance profiles are used to compare different methods (see \cite{zitp,zitq,flg3,mat2,mat4,mat5} for more details of performance profiles). A performance profile is defined as the cumulative function $\rho(\tau)$ representing the performance ratio with respect to a given metric and a set of methods. Give a set of methods $\mathcal{SO}$ and a set of problems $\mathcal{P}$, let $\varsigma_{p,s}$ be the performance of solver $s$ on solving $p$. The performance ratio is defined as $r_{p,s}=\varsigma_{p,s}/\min_{s\in\mathcal{SO}}\varsigma_{p,s}$. The cumulative function $ \rho_s(\tau)~~(s\in\mathcal{SO})$ is defined as
 \bee
  \rho_s(\tau)=\frac{|\{p \in \mathcal{P}:r_{p,s}\leq \tau\}|}{|P|}.
  \eee
It is noticed from the numerical experiments that the performance metric is sensitive to the number and types of algorithms considered in comparison process (see \cite{mat2,mat4,mat5} and references mentioned there). So algorithms should be compared pairwise for computing performance profile. The output of a multi-objective optimization solver is a set of non-dominated points. To justify  how much well-distributed this set is, the following metrics are considered for computing performance profile.\\\\
\textbf{$\Delta$-spread metric:} Let $x^1,x^2,\dots,x^N$ be the set of points obtained by a solver $s$ for problem $p$ and let these points be sorted by $f_j(x^i)\leq f_j(x^{i+1})$ $(i=1,2,\dots,N-1)$. Suppose $x^0$ is the best known approximation of global minimum of $f_j$ and $x^{N+1}$ is the best known global maximum of $f_j$, computed over all the approximated Pareto fronts obtained by different solvers. Define $\bar{\delta_j}$ as the average of the distances $\delta_{i,j}$, $ i=1,2,\dots,N-1.$ For an algorithm  $s$ and a problem $p$, the spread metric $\Delta_{p,s}$ is
\bee
\Delta_{p,s}:=\underset{j\in\Lambda_m}{\max}\left(\frac{\delta_{0,j}+\delta_{N,j}+\Sigma_{i=1}^{N-1} |\delta_{i,j}-\bar{\delta_j}|}{\delta_{0,j}+\delta_{N,j}+(N-1)\bar{\delta_j}}\right).
\eee
\textbf{Hypervolume metric:} Hypervolume metric of an approximate Pareto front with respect to a reference point $P_{ref}$ is defined as the volume of the total region dominated by the efficient solutions obtained by a method with respect to the reference point. Since it is very difficult to compute exact hypervolume metric, several techniques for computing hypervolume indicator are introduced by several researchers. Here $10000$ uniformly distributed random points are generated between $P_{ref}$ and ideal vector. Hypervolume metric is defined as $hv_{p,s}=N_{dom}/10000$, where  $N_{dom}$ is the number of points dominated by approximate Pareto front. Higher values of $hv_{p,s}$ indicate better performance using hypervolume metric. So while using the performance profile of the solvers measured by hypervolume metric we need to set $\widetilde{hv_{p,s}}=\frac{1}{hv_{p,s}}$.\\

\begin{figure}[!htbp]
    \centering
    \begin{subfigure}[b]{.47\textwidth}
        \centering
        \includegraphics[height=2.75cm,width=\linewidth]{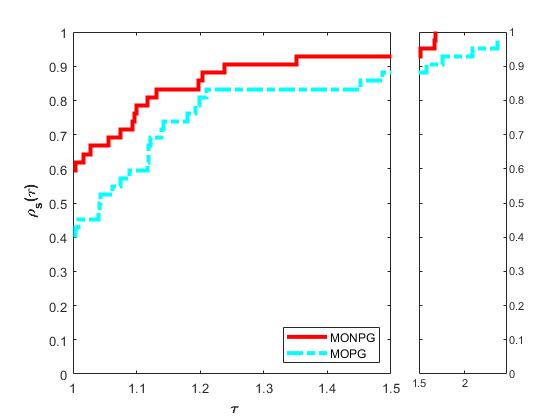}
        \caption{Between MONPG and MOPG}
        \label{fdel1}
    \end{subfigure}
    \begin{subfigure}[b]{.47\textwidth}
        \centering
        \includegraphics[height=2.75cm,width=\linewidth]{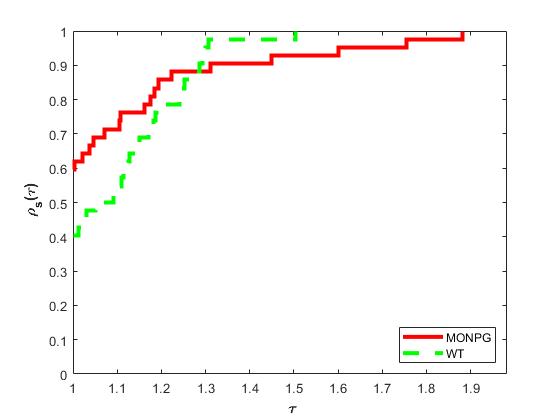}
        \caption{Between MONPG and WS}
        \label{fdel2}
    \end{subfigure}
\caption{Performance profiles using $\Delta$-spread metric}
\label{ctt1}
\end{figure}

\begin{figure}[!htbp]
    \centering
    \begin{subfigure}[b]{.47\textwidth}
        \centering
        \includegraphics[height=2.75cm,width=\linewidth]{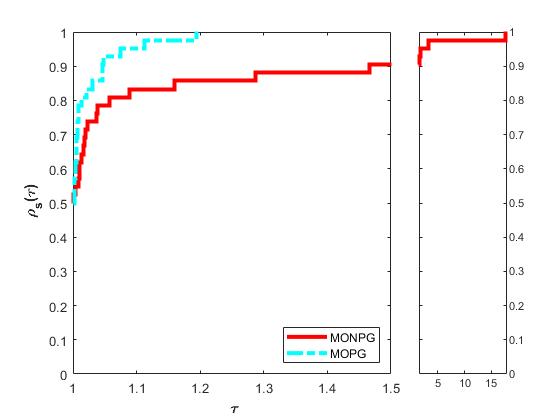}
        \caption{Between MONPG and MOPG}
        \label{fhv1}
    \end{subfigure}
    \begin{subfigure}[b]{.47\textwidth}
        \centering
        \includegraphics[height=2.75cm,width=\linewidth]{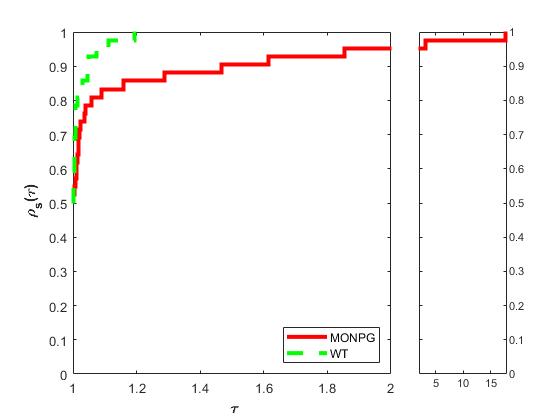}
        \caption{Between MONPG and WS}
        \label{fhv2}
    \end{subfigure}
\caption{Performance profiles using hypervolume metric}
\end{figure}

Performance profiles using $\Delta$-spread metric between MONPG and MOPG and MONPG and WS are provided in Figure \ref{fdel1} and \ref{fdel2} respectively. Performance profiles using hypervolume metric between MONPG and MOPG and MONPG and WS are provided in Figure \ref{fhv1} and \ref{fhv2} respectively.\\
Apart from these two metrics, performance profiles are computed using the number of iterations and total function evaluations. MONPG and WS use Hessian information of smooth function, whereas MOPG does not. Number of gradients (in each method) and Hessian (in MONPG and WS) evaluation is equal to the number of iterations. We have used forward difference formula to find the gradient and Hessian. This requires $n $function evaluations to find a gradient and  $\frac{1}{2}n(n+1)$ function evaluations to find a Hessian. So total number of function evaluations in MONPG and WS are
\bee
\# f=\#f+n\# it+\frac{1}{2}n(n+1) \# it
\eee
and number of function evaluations in MOPG is
\bee
\# f=\#f+n\# it
\eee
where $\# f$ and $\# it$ are number of function evaluations and the number of iterations respectively.\\

\begin{figure}[!htbp]
    \centering
    \begin{subfigure}[b]{.47\textwidth}
        \centering
        \includegraphics[height=2.75cm,width=\linewidth]{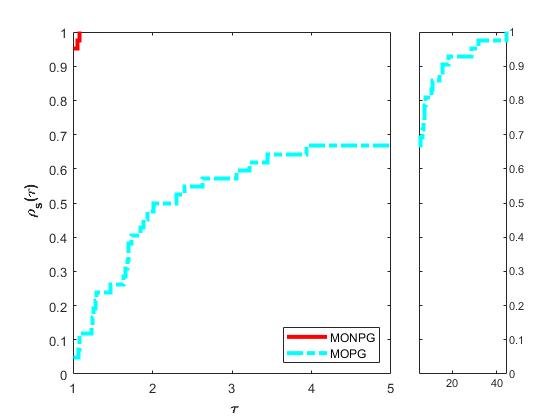}
        \caption{Between MONPG and MOPG}
        \label{fit1}
    \end{subfigure}
    \begin{subfigure}[b]{.47\textwidth}
        \centering
        \includegraphics[height=2.75cm,width=\linewidth]{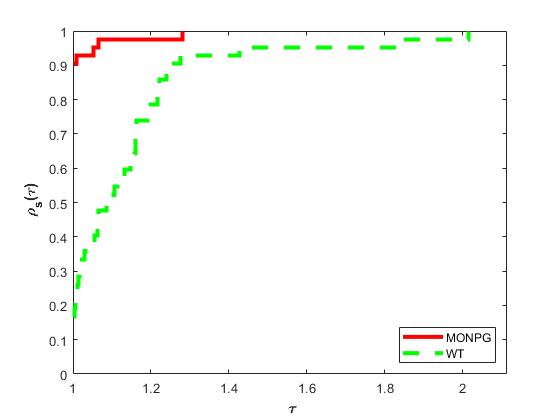}
        \caption{Between MONPG and WS}
        \label{fit2}
    \end{subfigure}
\caption{Performance profiles using number of iterations}
\end{figure}

\begin{figure}[!htbp]
    \centering
    \begin{subfigure}[b]{.47\textwidth}
        \centering
        \includegraphics[height=2.75cm,width=\linewidth]{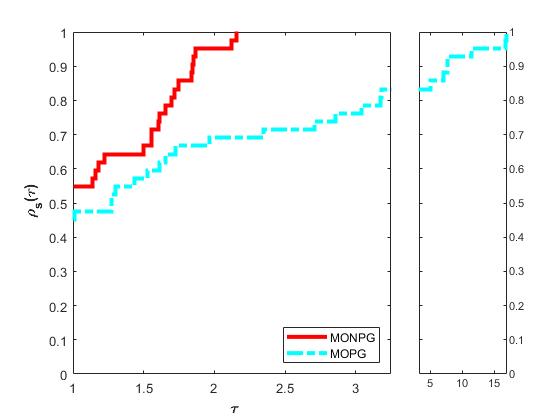}
        \caption{Between MONPG and MOPG}
        \label{ffun1}
    \end{subfigure}
    \begin{subfigure}[b]{.47\textwidth}
        \centering
        \includegraphics[height=2.75cm,width=\linewidth]{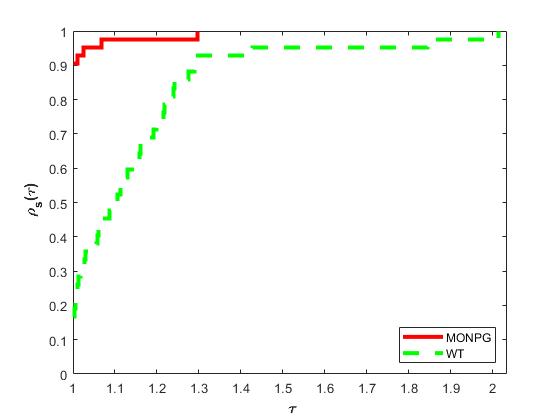}
        \caption{Between MONPG and WS}
        \label{ffun2}
    \end{subfigure}
\caption{Performance profiles using number of iterations}
\label{ctt2}
\end{figure}

Performance profiles using the number of iterations between MONPG and MOPG and MONPG and WS are provided in Figure \ref{fit1} and \ref{fit2} respectively. Performance profiles using hypervolume metric between MONPG and MOPG and MONPG and WS are provided in Figure \ref{ffun1} and \ref{ffun2} respectively.
\newpage
From performance profile figures (Figures \ref{ctt1}-\ref{ctt2}) one can observe that MONPG takes a less number of iterations and function evaluations in most cases and MONPG provides better results in $\Delta$-spread metric in most cases. 
\section{Conclusion}
In this paper, we have developed a Newton-type proximal gradient method for nonlinear convex multi-objective optimization problems. This method is free from any kind of priori chosen parameters or ordering information of objective function. Global convergence of the proposed method is justified under some mild assumptions. However, generating a well distributed approximate Pareto front is not addressed here. Generating a well distributed Pareto front of multi-objective optimization problems is not an easy task. One technique does not work for different type problems. In literature different techniques are used by several researchers such as initial point selection technique (\cite{flg3,mat5}), continuation method (\cite{sdd1,hillermeier1,martinb1}), and predictor corrector method (\cite{martin1,martin2}) etc. These methods are restricted to smooth multi-objective optimization problems. The extension of these ideas to composite multi-objective optimization problems is left as a scope of future research. Some hybrid methods can be developed by combining proximal gradient methods with genetic algorithms to generate well distributed approximate Pareto fronts, which is left as scope of future research. Apart from spreading of Pareto fronts, the proposed method is restricted to convex multi-objective optimization problems and Hessian of each smooth function is used in every iteration which is computationally expensive. In future we want to use the ideas of quasi-Newton methods to overcome these limitations.
\section*{Acknowledgements}
The author would like to thank the referees for their detailed comments and suggestions that have significantly improved the content as well as the presentation of the results in the paper.

\end{document}